\documentstyle{article}
\newtheorem{prethm}{{\bf Theorem}}

\newenvironment{thm}{\begin{prethm}{\hspace{-0.5
               em}{\bf.}}}{\end{prethm}}
\newtheorem{prelemma}{{\bf Lemma}}

\newenvironment{lemma}{\begin{prelemma}{\hspace{-0.5
               em}{\bf.}}}{\end{prelemma}}
\newtheorem{preex}{{\bf Example}}

\newenvironment{ex}{\begin{preex}{\hspace{-0.5
               em}{\bf.}}}{\end{preex}}
\newtheorem{preexmpl}{{\bf Example}}

\newenvironment{exmpl}{\begin{preexmpl}{\hspace{-0.8
               em}{\bf.}}}{\end{preexmpl}}
\newtheorem{prepro}{{\bf Proposition}}

\newenvironment{pro}{\begin{prepro}{\hspace{-0.5
               em}{\bf.}}}{\end{prepro}}
\newtheorem{precor}{{\bf Corollary}}

\newenvironment{cor}{\begin{precor}{\hspace{-1
               em}{\bf.}}}{\end{precor}}
\newtheorem{precorr}{{\bf Corollary}}

\newenvironment{corr}{\begin{precorr}{\hspace{-0.5
               em}{\bf.}}}{\end{precorr}}
\newtheorem{predeff}{{\bf Definition}}

\newenvironment{deff}{\begin{predeff}{\hspace{-0.85
               em}{\bf.}}}{\end{predeff}}
\newtheorem{preque}{{\bf Question}}

\newenvironment{que}{\begin{preque}{\hspace{-0.85
               em}{\bf.}}}{\end{preque}}
\newtheorem{preconj}{{\bf Conjecture}}

\newtheorem{preremark}{{\bf Remark}}

\newtheorem{prerem}{{\bf Remark}}

\newenvironment{rem}[1]{\begin{prerem}
               {\hspace{-0.7em}{\bf.}\hspace{-0.4em}}
               {\rm #1}}{\end{prerem}}
\newtheorem{preprob}{{\bf Problem}}

\newtheorem{preprobb}{{\bf Problem}}

\newenvironment{probb}{\begin{preprobb}{\hspace{-0.5
               em}{\bf.}}}{\end{preprobb}}
\newtheorem{prelem}{{\bf Theorem}}

\newenvironment{lem}{\begin{prelem}{\hspace{-0.5
               em}{\bf.}}}{\end{prelem}}
\newtheorem{preproof}{{\bf Proof.}}

\newenvironment{proof}[1]{\begin{preproof}{\rm
               #1}\hfill{$\Box$}}{\end{preproof}}

\newtheorem{presolution}{{\bf Solution.}}

\def\emline#1#2#3#4#5#6{%
       \put(#1,#2){\special{em:moveto}}%
       \put(#4,#5){\special{em:lineto}}}
\def\newpic#1{}

\newcommand{\utwolc}{{\rm U$2$LC}}
\newcommand{\uthrlc}{{\rm U$3$LC}}
\newcommand{\uklc}{{\rm U$k$LC}}
\newcommand{\m}{\mbox{\rm m}\,}
\newcommand{\mnum}{{\rm m}--number}
\def\Chi{\lower-.3ex\hbox{$\large\chi$}}
\title{\LARGE\sf On Uniquely List Colorable Graphs\footnote{
The research is partially supported by the
 Institute for Studies in Theoretical Physics and Mathematics (IPM),
Tehran, Iran.}}

\author{
{\rm M.\,Ghebleh} and {\rm E.S.\,Mahmoodian}\\
\vspace{2mm} \\
Institute for studies in theoretical Physics\\and Mathematics (IPM)\\
and\\
Department of Mathematical Sciences\\
Sharif University of Technology\\
P.O.\,Box 11365--9415\\
Tehran,\ \ I.\,R.\,Iran
}
\date{}
\begin{document}
\maketitle
\addtolength{\baselineskip}{1mm}
\begin{abstract}
Let  $G$ be a graph with $n$ vertices and suppose that for each vertex
$v$ in $G$, there exists a list of $k$ colors, $L(v)$, such that there
is a unique proper coloring for $G$ from this collection of lists,
then $G$ is called a {\it uniquely $k$--list colorable graph}.
Recently M.~Mahdian and E.S.~Mahmoodian characterized uniquely
$2$--list colorable graphs. Here we state some results which will pave
the way in characterization of uniquely $k$--list colorable graphs.
There is a relationship between  this concept and defining sets
in graph colorings and critical sets in latin squares.\\[.5cm]
\end{abstract}
\section{Introduction and preliminaries} 

We consider simple graphs which are finite, undirected,
with no loops or
multiple edges. For the necessary definitions and
notation we refer the reader to standard
texts, such as  \cite{west}.
In this section we mention some of the definitions and results
which are referred to throughout the paper.

For each vertex $v$ in a graph $G$, let $L(v)$ denote a list of colors
available for~$v$. A {\sf list coloring} from the given collection
of lists is
a proper coloring~$c$ such that $c(v)$ is chosen from $L(v)$.
We will refer to such a coloring as an {\sf $L$--coloring}.
The idea of list colorings of graphs is due independently to
V.~G.~Vizing~\cite{vizing} and to
P.~Erd\"{o}s, A.~L.~Rubin, and H.~Taylor~\cite{erdos}.
For a recent survey on list coloring we refer the interested reader
to N.~Alon~\cite{Alon}.
It is interesting to note that a list coloring of $K_n$ is nothing but a
system of distinct representatives (SDR) for the collection
${\cal L}=\{L(v)|v\in V(K_n)\}$.

Let $G$ be a graph with $n$ vertices and suppose that for each vertex
$v$ in $G$, there exists a list of $k$ colors $L(v)$, such that there
exists a unique $L$--coloring for $G$,
then $G$ is called a {\sf uniquely $k$--list colorable graph} or a
{\sf \uklc\ graph} for short.
\begin{exmpl}
\label{k4-e}
The graph $K_4\setminus e$ is a uniquely $2$--list colorable graph.
\end{exmpl}
In Figure~\ref{fk4-e}
a collection of lists is given, each of size two, and it can easily
be checked that there is a unique coloring with these lists.
\begin{figure}[th]
\begin{center}
\unitlength=.50mm
\special{em:linewidth 0.4pt}
\linethickness{0.4pt}
\begin{picture}(140.50,75.00)(10,20)
\put(104.00,52.00){\circle*{4.00}}
\put(79.00,27.00){\circle*{4.00}}
\put(54.00,52.00){\circle*{4.00}}
\put(79.00,77.00){\circle*{4.00}}
\
\put(110.00,52.00){\makebox(0,0)[cc]{$13$}}
\put(79.00,21.00){\makebox(0,0)[cc]{$12$}}
\put(48.00,52.00){\makebox(0,0)[cc]{$12$}}
\put(79.00,83.00){\makebox(0,0)[cc]{$13$}}
\
\emline{104.00}{52.00}{1}{79.00}{27.00}{2}
\emline{104.00}{52.00}{1}{54.00}{52.00}{2}
\emline{104.00}{52.00}{1}{79.00}{77.00}{2}
\emline{79.00}{27.00}{1}{104.00}{52.00}{2}
\emline{79.00}{27.00}{1}{54.00}{52.00}{2}
\emline{54.00}{52.00}{1}{104.00}{52.00}{2}
\emline{54.00}{52.00}{1}{79.00}{27.00}{2}
\emline{54.00}{52.00}{1}{79.00}{77.00}{2}
\emline{79.00}{77.00}{1}{104.00}{52.00}{2}
\emline{79.00}{77.00}{1}{54.00}{52.00}{2}
\end{picture}
\caption{$K_4\setminus e$}
\end{center}
\label{fk4-e}
\end{figure}
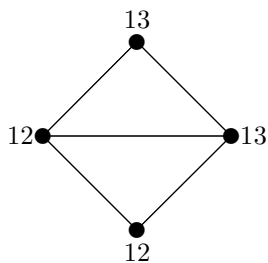
\begin{rem}{
\hspace{-.2ex}It is clear from the definition of uniquely $k$--list
colorable graphs that each \uklc\ graph is also a U$(k-1)$LC graph.
}\end{rem}

The following theorem of Marshal~Hall, which
is a corollary of the celebrated Marriage Theorem of P.~Hall and gives
a lower bound for the number of SDRs,  is a motivation for the
definition of \uklc\ graphs.
\begin{lem}~{\rm\cite{hall}}
\label{mhall}
If $n$ sets $S_1, S_2, \ldots,S_n$
have an SDR and the smallest of
these sets contains $k$ objects, then if $k\ge n$, there are at least
$k(k-1) \cdots (k- n +1)$ different SDRs; and if $k< n$, there
are at least $k!$ different SDRs.
\end{lem}
\begin{cor} \label{csds}
If the sets $S_1, S_2, \ldots,S_n$ have an
SDR and the smallest of these sets is of
size $k$ {\rm (} $k > 1${\rm ),} then
they have at least two  SDRs. Or equivalently,
the complete graph $K_n$ is not \uklc.
\end{cor}
If in the above corollary instead of $K_n$ we take any graph, then
it is natural to ask the following question.
\begin{que}
\label{forwhich}
For which graphs does the result of the above corollary hold?
\end{que}
We say that a graph $G$
has {\sf the property $M(k)$} ($M$ for Marshal~Hall)
if and only if it is not uniquely $k$--list colorable. So $G$ has
the property $M(k)$ if for any collection of lists
assigned to its vertices, each of size $k$,
either there is no list coloring for $G$ or
there exist at least two list colorings. Note that if one tries to
relate the idea of uniqueness to list coloring, then he or she reaches
this definition naturally.

M.~Mahdian and E.S.~Mahmoodian characterized uniquely
$2$--list colorable graphs. They showed that,
\begin{lem}~{\rm\cite{mahmah}}
\label{u2lc}
A connected  graph $G$ has the property $M(2)$ if and only if
every block of $G$ is either
a cycle, a complete graph, or a complete bipartite graph.
\end{lem}
\noindent
It seems that characterizing \uklc\ graphs for any $k$
is not that easy. Even the \uthrlc\ graphs
seem to be  difficult to characterize.
For example it will be shown below that, while there are some complete
tripartite graphs which have the property $M(3)$,
the property does not hold for any complete tripartite graph.

The following definition was first given in~\cite{tabriz}.
\begin{deff}
The {\sf \mnum\ } of a graph $G$, denoted by $\m(G)$, is defined
to be the least integer $k$ such that $G$ has the property $M(k)$.
\end{deff}
E.\,S.~Mahmoodian and M.~Mahdian in~\cite{tabriz} have obtained some
results on the \mnum\ of planar graphs and introduced some upper bounds
on~$\m(G)$.

It is obvious from the definition of a \uklc\ graph that the graph $G$ is
\uklc\ if and only if $k<\m(G)$.
For example, one can easily see that the graph $K_4\setminus e$ has the
property $M(3)$ and in the above example we saw that it is \utwolc, so
$m(K_4\setminus e)=3$.

The concept of \uklc\ graphs also arise naturally in finding
defining sets of graphs.
In a given graph $G$, a set of vertices $S$ with an assignment of colors
is called a {\sf defining set of $k$--coloring}, if there exists a
unique extension of the colors of $S$ to a $ k $--coloring
of the vertices of $G$.
For more information on defining sets see~\cite{mnz}. As it is
mentioned there, critical sets in latin squares are just the minimal
defining sets of $n$--colorings of $K_n\times K_n$.
A {\sf latin square} is an $n \times n$ array from
the numbers $1, 2, \ldots, n$ such that each of these numbers occurs
in each row and in each column exactly once. A {\sf critical set }
in an $n \times n$ array is a set $S$ of  entries, such that
there exists a unique extension of  $S$ to  a latin square of size
$n$ and no proper subset of $S$ has this property. For a survey on
critical sets in latin squares see~\cite{keedwell}.

Each set of vertices $S\subset V(G)$ with an assignment of colors
induces a list of colors for each vertex in
$G\setminus S$. So to find out if $S$ is a defining set or not,
we need to know whether $G\setminus S$
is uniquely list colorable with those lists.

In this paper we state some results which are towards characterizing
\uklc\ graphs. In Section 2 we introduce some results which are helpful
in  determining the \mnum\ of some graphs. In Section 3 some theorems about
complete multipartite graphs are discussed. In Section 4 we present some
examples of \uklc\ graphs, and finally in the last section we pose some
open problems.

\section{Some general results} 
The following lemma is very useful throughout the paper.
\begin{lemma}
For every graph $G$ we have $\m(G)\le |E(\overline{G})|+2$.
\label{r5}
\end{lemma}
\begin{proof}{
Proof is by induction on $r=|E(\overline{G})|$. In the case $r=0$,
$G$ is a complete graph and by Theorem~\ref{u2lc} it has the
property $M(2)$. Assume that the statement is true for every graph
$H$ with $|E(\overline{H})|<r$ and let $G$ be a graph whose
complement has $r$ edges. Suppose that there are assigned some
lists of colors $L(w)$ of size at least $r+2$ to the vertices of
$G$ and $G$ has an $L$--coloring $c$. Let $u$ and $v$ be two
nonadjacent vertices of $G$. To obtain another $L$--coloring for
$G$, we consider two cases.

If $c(u)\not=c(v)$, consider the graph $G_1=G+uv$. We have
$|E(\overline{G_1})|=r-1$, and by induction hypothesis $\m(G_1)\le r+1$.
So there exists another $L$--coloring for $G_1$ which is also legal for
$G$ itself.

Now if $c(u)=c(v)$, consider the graph $G_2=G\setminus\{w|c(w)=c(v)\}$.
If $V(G_2)=\emptyset$, then $G$ is a null graph and the statement is
trivial.
Otherwise we have $|E(\overline{G_2})|<r$, and by induction hypothesis
$\m(G_2)\le r+1$. Assign to each vertex $w$ of $G_2$ the list
$L'(w)=L(w)\setminus\{c(u)\}$. Again $c|_{V(G_2)}$ is an
$L'$--coloring for $G_2$ and since $|L'(w)|\ge r+1$ for each $w\in
V(G_2)$, there exists another $L'$--coloring for $G_2$ which can
be extended to an $L$--coloring of $G$, different from $c$, by
giving the color $c(u)$ to all the vertices of $G$ which are not
in $G_2$.%
}\end{proof}
From the following theorem we can deduce a lower bound for the number of
vertices in a \uklc\ graph.
\begin{thm}
\label{r6}
If a graph $G$ has at most $3k$ vertices, then $m(G)\le k+1$.
\end{thm}
\begin{proof}{
Proof is by induction on $k$. For $k=1$ the statement obviously holds.
Suppose that $k\ge 2$, and $G$ is a graph with at most $3k$ vertices,
and let there be lists of colors, each
of size at least $k+1$, assigned to the vertices of $G$
and further suppose that there exists
a list coloring $c$ for $G$, from these lists. We show that there exists
another coloring for $G$ from these lists.

If one color class has at least three vertices, we can remove that
class from $G$ and its color from the lists of remaining vertices,
and by induction hypothesis a new coloring exists for the
remaining graph which extends to all of $G$. So assume that each
color class has at most two vertices. By adding new edges between
all vertices with different colors in $c$, we obtain a graph whose
complement is  union of some  $K_1$s and some $K_2$s. Denote the
number of $K_2$s by $r$. If $r\le k-1$, we obtain a new coloring
by the lemma above, otherwise $r\ge k$. Now if there exists a vertex
$v$ whose list contains a color $x$ which is not used in the
coloring $c$, then we can obtain a new coloring by changing the
color of $v$ to $x$. Otherwise the union of all lists has exactly
$n-r\le 2k$ elements. If $u$ and $v$ are two vertices such that
$c(u)=c(v)$, then since the unused colors in the lists of $u$ and
$v$ are chosen from a ($2k-1$)--set, thus $u$ and $v$ must have a
common unused color. Consider a $K_{n-r}$ obtained by identifying
all the vertices in each color class of $c$ to a vertex. The list
of each vertex in this $K_{n-r}$ is the intersection of the lists
of the vertices in the corresponding color class. So each list of
the vertices in $K_{n-r}$ has at least $2$ elements, and there
exists a coloring for it from these lists. Hence by the property
$M(2)$ of $K_{n-r}$ we obtain a new coloring on it, which gives a
new coloring for $G$. %
}\end{proof}
The following two corollaries are immediate from the theorem above.
The first one gives an upper bound for the \mnum\ of a graph and the
second one introduces a lower bound for the number of vertices in a
\uklc\ graph.
\begin{corr}
\label{r8}
If a graph $G$ has $n$ vertices then $m(G)\le\lceil n/3\rceil+1$.
\end{corr}
\begin{corr}
\label{r7}
Every \uklc\ graph has at least $3k-2$ vertices.
\end{corr}
Corollary~\ref{r7} implies that a necessary condition to have equality in
Lemma~\ref{r5} is  $|V(G)|\ge 3|E(\overline{G})|+1$. In the following
proposition we see that when the edges of $\overline{G}$ are independent
this condition is also sufficient .
\begin{pro}
\label{eqKnM}
If $F$ is a set of $r$ independent edges in $K_n$ and $n\ge 3r+1$,
then $m(K_n\setminus F)=r+2$.
\end{pro}
\begin{proof}{
Suppose $F=\{x_1y_1,\ldots,x_ry_r\}$, and $z_0,\ldots,z_s$
are the vertices in $K_n\setminus V(F)$. By the hypothesis $s\ge r$.
Assign the list $\{0,1,\ldots,r\}$ to each $x_i$ and to $z_0$, and
for each $i\ge 1$ assign the list $\{r+1,\ldots,2r,i\}$ to $y_i$, and
the list $\{1,\ldots,r,r+i\}$ to $z_i$.
Since the induced subgraph of $K_n\setminus F$ on $\{x_1,\ldots,x_r,z_0\}$
is a complete graph, all the colors $0,1,\ldots,r$ must appear on these
vertices in any coloring of $K_n\setminus F$ from the assigned lists.
So for each $i\ge 1$ the vertex $z_i$ must take the color $r+i$,
and for each $i\ge 1$ $y_i$ receives the color $i$. Finally each $x_i$
must take the color $i$, and $z_0$ takes the color $0$.
}\end{proof}
\section{Complete multipartite graphs}

It is shown in \cite{mahmah} that any complete bipartite graph has
the property $M(2)$. In the following theorem it is shown that one
can not expect similar statement for complete tripartite graphs.
\begin{thm}
For each $k\ge 2$, there exists a complete tripartite \uklc\ graph.
\end{thm}
\begin{proof}{Let $A=\{a_1,\dots,a_{k-1}\}$, $B=\{b_1,\dots,b_{k-1}\}$,
and $C=\{c_1,\dots,c_{k-1}\}$ be mutually disjoint sets. We denote all
$(k-1)$-subsets of $B\cup C$ by $\{A_1,\dots,A_m\}$, those of $A\cup C$
by $\{B_1,\dots,B_m\}$, and those of $A\cup B$ by $\{C_1,\dots,C_m\}$;
where $m={2k-2 \choose k-1}$.

Now consider a complete tripartite graph $K_{m(k-1), m(k-1),m(k-1)}$
with the following list of colors on vertices in three parts,
respectively:
$A_i \cup \{a_j \}$,
$B_i \cup \{b_j \}$, and
$C_i \cup \{c_j \}$; where
$i= 1, 2, \dots, m$ and
$j= 1, 2, \dots, k-1$.
We show that there is a unique coloring for this graph from the assigned
lists.

First note that the union of all lists is $A \cup B \cup C$ which has
$3(k-1)$ elements.
We show that in any coloring of this graph, there are at least $k-1$ colors
present on the vertices of each part.
To show this, suppose to the contrary that there exists a coloring
in which one part uses less than $k-1$ colors. Without loss of generality
let $L$ be the set of colors used to color the first part, and $|L|<k-1$.
Then $(B \cup C) \setminus L$ has at least $k$ elements and
$A  \setminus L$ has at least one element.
Now consider a set $L'$ which contains $k-1$ elements from the set
$(B \cup C) \setminus L$  and an element from
$A  \setminus L$. Then $L \cap L' = \emptyset $. But there is a vertex
in the first part whose list is $L'$, a contradiction. So each part has
at least $k-1$ colors and since we have $3(k-1)$ colors altogether, thus
in any coloring each part has exactly $k-1$ colors. It can be easily
verified that the colors of each of the three parts must be $A$, $B$,
and $C$, respectively. Therefore there is a unique coloring for
$K_{m(k-1),m(k-1),m(k-1)}$ from the assigned lists.
}\end{proof}
The following theorem and the propositions which follow are preparations
to prove our main theorem of this section, Theorem \ref{charu3lc}, which is a
characterization of uniquely $3$--list colorable complete multipartite
graphs except for finitely many of them. The proof of the following useful
lemma is immediate.
\begin{lemma}
If $L$ is a $k$--list assignment to the vertices in the graph $G$, and $G$
has a unique $L$--coloring, then $|\bigcup_vL(v)|\ge k+1$ and all these
colors are used in the (unique) $L$--coloring of $G$.
\end{lemma}
\begin{thm}
\label{r12}
If $G$ is a complete multipartite graph which has an induced \uklc\
subgraph, then $G$ is \uklc.
\end{thm}
\begin{proof}{
Let $H$ be an induced subgraph of $G$ which is \uklc. Assume that
$L$ is a $k$--list assignment to the vertices in $H$, by which $H$ has a
unique list coloring. For the vertices in $G$ we introduce
lists of colors each of size $k$, such that $G$ is uniquely colorable
by these lists.
Assign the list $L(v)$ to each vertex $v$ in $H$. For each part of $G$ that
contains some vertices in $H$, consider a vertex $v$ in $H$ in that part
and assign the list $L(v)$ to all vertices in $G\setminus V(H)$ in that
part.
In any part of $G$ which does not contain any vertex in $H$, we assign a
list $A\cup \{i\}$, where $A$ is a set of $k-1$ colors from the $L$--coloring
of $H$ and $i$ is a new color.
}\end{proof}
We use the notation $K_{s*r}$ for a complete $r$--partite graph in which
each part is of size $s$. Notations such as $K_{s*r,t}$, etc. are used
similarly.
\begin{pro}
\label{kPartU3LC}
The graphs
$K_{3,3,3}$, $K_{2,4,4}$, $K_{2,3,5}$, $K_{2,2,9}$,
$K_{1,2,2,2}$, $K_{1,1,2,3}$, $K_{1,1,1,2,2}$,
$K_{1*4,6}$, $K_{1*5,5}$, and $K_{1*6,4}$
are \uthrlc.
\end{pro}
\begin{proof}{
First we show  the truth of the statement for $K_{1,1,2,3}$ and $K_{1*4,6}$.

For $K_{1,1,2,3}$, let $\{a\}$, $\{b\}$, $\{c,d\}$, and $\{e,g,f\}$ be
the parts in $K_{1,1,2,3}$.
We assign the following lists for the vertices of this graph:
$L(a)=L(c)=L(f)=\{1,3,4\}$, $L(b)=L(d)=L(g)=\{2,3,4\}$,
and $L(e)=\{1,2,4\}$.
A unique coloring exists from the assigned lists, because
the vertices $b$, $d$ and $g$ form a triangle and all of them have the
list $\{2,3,4\}$, thus the colors $2$, $3$, and $4$ all occur on
these vertices.
The vertex $a$ is adjacent to these three vertices, so it is forced to take
the color $1$. Now the colors $3$, $4$ must both occur on $c$ and $f$
so $b$ must take the color $2$.
Finally $e$ is forced to take the color $4$, $c$ and $d$ must take $3$,
and the two remaining vertices $f$ and $g$ must take $4$.

For $K_{1*4,6}$, assign the lists $\{1,5,6\}$, $\{2,5,6\}$,
$\{3,5,6\}$, and $\{4,5,6\}$ to the vertices in the parts which
have one vertex each, and the lists
$\{1,2,5\}$, $\{1,3,5\}$, $\{1,4,6\}$, $\{2,3,6\}$, $\{2,4,6\}$,
$\{3,4,5\}$ to the vertices in the last part. In any coloring we
need all six colors because the last part needs at least two colors.
Now none of the colors $1, 2, 3$, and $4$ can appear on the last part
because in that case we need more than two colors on the last part,
a contradiction.

For each of the other eight graphs one can check by similar argument
that it has a unique coloring from the lists given below:

\noindent
$K_{3,3,3}$:
$\{\sf
\{{\underline 1}34, {\underline 1}35, {\underline 2}45\},
\{12{\underline 3}, 1{\underline 4}5, {\underline 3}56\},
\{13{\underline 6}, 14{\underline 5}, 23{\underline 5}\}
\}$\\
$K_{2,4,4}$:
$\{\sf
\{{\underline 1}35, {\underline 2}46\},
\{1{\underline 3}5, 2{\underline 4}6, {\underline 3}56, {\underline 4}56\},
\{12{\underline 5}, 34{\underline 5}, 14{\underline 6}, 23{\underline 6}\}
\}$\\
$K_{2,3,5}$:
$\{\sf
\{{\underline 1}46, {\underline 2}35\},
\{1{\underline 3}6, 2{\underline 3}5, {\underline 4}56\},
\{12{\underline 5}, 34{\underline 5}, 13{\underline 6}, 23{\underline 6},
24{\underline 6}\}
\}$\\
$K_{2,2,9}$:
$\{\sf
\{{\underline 1}56, {\underline 2}34\},
\{1{\underline 3}5, 1{\underline 4}6\},
\{12{\underline 5}, 13{\underline 5}, 14{\underline 5}, 12{\underline 6},
13{\underline 6}, 14{\underline 6},
24{\underline 5}, 34{\underline 5}, 23{\underline 6}\}
\}$\\
$K_{1,2,2,2}$:
$\{\sf
\{{\underline 1}23\},
\{1{\underline 2}3, {\underline 2}45\},
\{12{\underline 3}, {\underline 3}45\},
\{12{\underline 4}, 12{\underline 5}\}
\}$\\
$K_{1,1,1,2,2}$:
$\{\sf
\{{\underline 1}45\},
\{{\underline 2}45\},
\{{\underline 3}45\},
\{12{\underline 4}, 3{\underline 4}5\},
\{12{\underline 5}, 34{\underline 5}\}
\}$\\
$K_{1*5,5}$:
$\{\sf
\{{\underline 1}67\},
\{{\underline 2}67\},
\{{\underline 3}67\},
\{{\underline 4}67\},
\{{\underline 5}67\},
\{12{\underline 6}, 34{\underline 6}, 15{\underline 6}, 25{\underline 7},
34{\underline 7}\}
\}$\\
$K_{1*6,4}$:
$\{\sf
\{{\underline 1}78\},
\{{\underline 2}78\},
\{{\underline 3}78\},
\{{\underline 4}78\},
\{{\underline 5}78\},
\{{\underline 6}78\},
\{12{\underline 7}, 34{\underline 7}, 12{\underline 8}, 56{\underline 8}\}
\}$
}\end{proof}
\begin{pro}
\label{r10}
$m(K_{2,2,3})=m(K_{2,3,3})=3$.
\end{pro}
\begin{proof}{
By Theorem~\ref{u2lc} the graph $K_{2,2,3}$ is a \utwolc\ graph,
so $m(K_{2,2,3})\ge 3$. We show that $m(K_{2,2,3})=3$. Suppose
that there are assigned color lists, each of size at least $3$, to
the vertices in $K_{2,2,3}$ and $c$ is a coloring from those
lists. If all vertices in a part of $K_{2,2,3}$ have the same
color in $c$, we can remove that color from the lists of the other
two parts and by the property $M(2)$ of complete bipartite graphs
we obtain a different coloring on those parts which is extendible
to $K_{2,2,3}$. So suppose that at least two colors appear on each
part. Add new edges between those nonadjacent vertices that take
different colors in $c$, the resulting graph is a $K_7$ or
$K_7\setminus e$, both of which have the property $M(3)$. So we
obtain another coloring which is a legal coloring for $K_{2,2,3}$.

The second graph is checked by a computer program and it has the property
$M(3)$, so by Theorem~\ref{u2lc} its \mnum\ is equal to 3.
}\end{proof}
\begin{pro}
\label{k1rs}
 Every complete tripartite graph $K_{1,s,t}$ has the property $M(3)$. Thus if
\ \ $\max \{s, t\} \geq 2 $, then \ \ $m(K_{1,s,t}) = 3$.
\end{pro}
\begin{proof}{
The proof is immediate by a technique similar to one used in
Proposition~\ref{r10}.
}\end{proof}
\begin{pro}
\label{K111r}
For each $s\ge 2$, $m(K_{1,1,1,s})=3$.
\end{pro}
\begin{proof}{
Suppose for each $v\in V(K_{1,1,1,s})$ there is assigned a color
list $L(v)$ of size $3$, and $K_{1,1,1,s}$ has an $L$--coloring
$c$. If one of the vertices in $K_{1,1,1,s}$ has a color in its
list which is not used in $c$, we obtain a new $L$--coloring for
$K_{1,1,1,s}$ by simply putting that unused color on that vertex.
So suppose that each color in $\cup_v L(v)$ is used in the
coloring.

Call the vertices in the first three parts $x, y$, and $z$, and the
vertices in the last part $w_1,\ldots,w_s$. Suppose that the colors
of $x, y$, and $z$ in the coloring are $1, 2$, and $3$, respectively.
So for each $i$, $L(w_i)$ contains $c(w_i)$ and two colors from
$1, 2$, and $3$.

If two of the vertices $x, y,$ and $z$, say $x$ and $y$ have some
colors of the last part in their lists, $c(w_p)\in L(x)$ and
$c(w_q)\in L(y)$ where $c(w_p)\not=c(w_q)$, then we obtain a new
coloring $c'$ for $K_{1,1,1,s}$ by putting $c(w_p)$ on $x$, $c(w_q)$ on $y$,
$c(z)$ on $z$, and since for each $i=1,2,\ldots,s$, there exists
$c'(w_i)\in L(w_i)\cap\{1,2\}$, we change each $c(w_i)$ by this
$c'(w_i)$.
Otherwise, either there is at most one color of the last part in
$L(x)\cup L(y)\cup L(z)$, or there is one of $x, y,$ and $z$,
say $x$, whose list contains two colors from the last part, and
two other have no color of the last part in their lists.
In the former case we can obtain a new coloring for the triangle
induced on $x,y,$ and $z$ from the lists $L(v)\cap\{1,2,3\}$ on
each $v\in\{x,y,z\}$, by the property $M(2)$ of $K_3$.
In the latter case a new coloring can be obtained by replacing the
colors of $y$ and $z$.

We showed that $K_{1,1,1,s}$ has the property $M(3)$, and so $m(K_{1,1,1,s})\le 3$. On the
other hand it has an induced $K_{1,1,2}$ subgraph which is a \utwolc\
graph, and so we have $m(K_{1,1,1,s})>2$.
}\end{proof}
\begin{pro}
\label{K1s3}
For every $r\ge 2$, we have $m(K_{1*r,3})=3$.
\end{pro}
\begin{proof}{
Suppose there are some lists of colors each of size $3$ assigned to the
vertices of $K_{1*r,3}$, which have a coloring.
We consider two cases and in each case obtain a new coloring for $K_{1*r,3}$
from these lists.
First consider the case that all vertices in the last part
take the same color in the given coloring. By removing this color
from the lists of other vertices, they have a new coloring because
the complete graphs have the property $M(2)$. So at least
two colors appear on the vertices in last part. Add new edges between
those vertices in the last part that have different colors. The resulting
graph is either a complete graph or a complete graph with an edge removed,
and we know that both of those graphs have the property $M(3)$.
So a new coloring can be obtained from the lists for the new graph.
This coloring is also valid for $K_{1*r,3}$.
}\end{proof}
Now we state our main theorem of this section.
\begin{thm}
\label{charu3lc}
Let $G$ be a complete multipartite graph that is not $K_{2,2,r}$,
for $r=4,5,\ldots,8$, $K_{2,3,4}$, $K_{1*4,4}$, $K_{1*4,5}$,
or $K_{1*5,4}$ then $G$ is \uthrlc\ if and only if it has one of the
graphs in Proposition~{\rm\ref{kPartU3LC}} as an induced subgraph.
\end{thm}
\begin{proof}{
If $G$ has one of the
graphs of Proposition~{\rm\ref{kPartU3LC}} as an induced subgraph, then
it is \uthrlc\ by Theorem~\ref{r12}.
So we prove the other side of the statement. Assume that $G$
is not one of the graphs mentioned in the statement and it does not have
any graphs of Proposition~{\rm\ref{kPartU3LC}} as an induced subgraph.
We show that it is not \uthrlc. There are two cases to be considered.

\noindent
(i)  $G = K_{1*r, s}$, for some $r$ and $s$.
  If $r \leq 3$ or $s \leq 3$, then by Proposition~\ref{K111r} and
  Proposition~\ref{K1s3} it has the property $M(3)$. So assume
  $r \geq 4$ and $s \geq 4$. Since $G$ does not contain a $K_{1*4,6}$
  we must have $4 \leq s \leq 5$. If $s=5$ we have $r=4$ which is exempted.
  If $s=4$ we have $r=4$ or 5, which are also exempted.
%

\noindent
(ii)  $G$ has at least two parts whose sizes are greater than 1. Since it
  does not contain a $K_{1,1,1,2,2}$, it is either 4--partite,
  tripartite, or bipartite.

\noindent
  If $G$ is bipartite, it is not \uthrlc, by Theorem~\ref{u2lc}.

\noindent
  If $G$ is 4--partite, since it does not contain a $K_{1,2,2,2}$ or a
  $K_{1,1,2,3}$, it must be  $K_{1,1,2,2}$ which is not \uthrlc\ by
  Theorem~\ref{r6}.

\noindent
  So assume that $G = K_{r,s,t}$ for some $t \le s \le r$.
  Since it does not contain a $K_{3,3,3}$ we have $t \le 2$.
  If $t=1$ then it is not \uthrlc\ by Proposition~\ref{k1rs}.

\noindent
  If $t=2$, since it does not
  contain a $K_{2,4,4}$ we must have $s \le 3$.

\noindent
  If $s=2$ then
  $G$ must be a $K_{2,2,r}$ with $r \le 8$. But now, if $r \le 3$ it is not
  \uthrlc\ by Proposition~\ref{r10}, and the cases of $4 \le r \le 8$
  are exempted.

\noindent
  If $s=3$ then $G= K_{2,3,r}$ where $r \le 4$.
  Then if $r \le 3$ it is not \uthrlc, by Proposition~\ref{r10}, and
  for $r=4$ it is exempted.
}\end{proof}
\section{Some examples of \uklc\ graphs}
In this section we introduce some examples of \uklc\ graphs.
\begin{ex}
\label{mg23}
The graph $K_{1*k,2*(k-1)}$ has \mnum\ equal to $k+1$.
\end{ex}
\begin{proof}{
This is the example given in~\cite{mahmah} as a \uklc\ graph.
It is a special case of graphs discussed in Proposition~\ref{eqKnM}.
}\end{proof}
\begin{ex}
\label{r2}
The graph $K_{1,2*(k-1),k-1}$ has
\mnum\ $k+1$.
\end{ex}
\begin{proof}{
From each of the first $k$ parts choose a vertex and assign to it the list
$\{1,\ldots,k\}$. To the other vertex in $i$--th part ($2\le i\le k$)
assign the list $\{k+1,\ldots,2k-1,i\}$.
Finally in the last part, assign the list $\{1,\ldots,k-1,k+j\}$ to
the $j$--th vertex in that part ($1\le j\le k-1$).
Since this graph has a subgraph $K_k$ which has the list
$\{1,\ldots,k\}$ on each of its vertices, by a similar argument
as in the proof of Proposition~\ref{kPartU3LC}, a unique coloring from
these lists for $K_{1,2*(k-1),k-1}$ can be obtained.
}\end{proof}
\begin{ex}
\label{r11}
The complete $(k+1)$--partite graph $K_{1,1,2,\ldots,k}$ is \uklc.
\end{ex}
\begin{proof}{
We use the colors from the set $A=\{1,2,\ldots,k+1\}$. Assign
the list $A\setminus \{k\}$ to the vertex in the first
part, and in the $(i+1)$--th part ($1\le i\le k$) assign the
list $A\setminus \{k-j+2\}$ to the $j$--th vertex $(1\le j\le i)$.
Since $\Chi(K_{1,1,2,\ldots,k})=k+1$, we need $k+1$ colors
to color this graph, so all of the colors must be used and
in each part we must have exactly one color. Hence the
vertices in the $(k+1)$--th part must all take the color $1$,
the vertices in the $k$--th part must all take the color $2$,
\ldots, the single vertex in the second part must take the
color $k$, and finally the single vertex in the first part is
forced to take the color $k+1$.
}\end{proof}
\begin{ex}
\label{r3}
The graph ${\cal U}_k$ constructed below has \mnum\ $k+1$:

Let the set $\{v_1,\ldots,v_{3k-2}\}$ be the set of vertices in
${\cal U}_k$. The edges in ${\cal U}_k$ are $v_iv_j$s
$(i\neq j)$ where:
\begin{itemize}
\item{$1\le i,j\le k$,}
\item{$1\le i\le k$ and $k+1\le j\le 2k-1$,}
\item{$k+1\le i\le 2k-1$ and $2k\le j\le 3k-2$,}
\item{$1\le i\le k-1$ and $2k\le j\le 3k-i-1$.}
\end{itemize}
\end{ex}
\begin{proof}{
Assign the list $\{1,\ldots,k\}$ to $v_1,\ldots,v_k$,
the list $\{1,\ldots,k-1,i\}$ to $v_i$ where $k+1\le i\le 2k-1$,
and the list $\{k+1,\ldots,2k-1,i-2k+1\}$ to $v_i$ where $2k\le i\le 3k-2$.
Again since there exists a $K_k$ in ${\cal U}_k$ induced on
$\{v_1,\ldots,v_k\}$ and with a similar argument as in the proof of
Proposition~\ref{kPartU3LC}, a unique coloring from these lists for
${\cal U}_k$ is obtained.
}\end{proof}
\begin{ex}
\label{Tk}
The graph ${\cal T}_k$ constructed below is \uklc\ for each $k\ge2$:
$$V(G)=\{a_1,\ldots,a_{k-1},b_1,\ldots,b_k,c_1,\ldots,c_{k-1},
d_1,\ldots,d_{2k-3}\},$$
and for edges,
\begin{itemize}
\item {Make a $K_{2k-1}$ on $a_i$s and $b_i$s},
\item {Join $b_i$s to $c_i$s and $c_i$s to $d_i$s},
\item {Join $a_i$ to $d_j$ for $1\le i\le k-1$ and $i\le j\le k-1$},
\item {Join $b_i$ to $d_j$ for $3\le i\le k$ and $k\le j\le k+i-3$}.
\end{itemize}
\end{ex}
\begin{proof}{
Assign some lists to the vertices in ${\cal T}_k$ as follows:
$L(a_i)=\{1,\ldots,k\}$,
$L(b_1)=\{k,\ldots,2k-1\}$,
$L(b_i)=\{i-1,k+1,\ldots,2k-1\}$ for $i>1$,
$L(c_i)=\{k+1,\ldots,2k-1,2k+i-1\}$, and
$L(d_i)=\{i+1,2k,\ldots,3k-2\}$.
It is easy to check that ${\cal T}_k$ has a unique coloring from
these lists.
}\end{proof}
\section{Some open problems}
The following problems arise naturally from the work.
\begin{probb}
Verify the property $M(3)$ for the graphs exempted in
Theorem~{\rm\ref{charu3lc}}, i.e. $K_{2,2,r}$ for $r=4,5,\ldots,8$,
$K_{2,3,4}$, $K_{1*4,4}$, $K_{1*4,5}$, and $K_{1*5,4}$.

\end{probb}
%
%
\begin{probb}
\label{prob3}
Characterize all graphs with \mnum\ $3$.
\end{probb}
\begin{probb}
What is the computational complexity of the property $M(3)$?
\end{probb}
\section*{Acknowledgement}
We thank Bashir~Sadjad who pointed out that in Lemma~\ref{r5}, the edges
of $\overline{G}$ are not necessarily supposed to be independent.
\newcommand{\noopsort}[1]{} \newcommand{\printfirst}[2]{#1}
  \newcommand{\singleletter}[1]{#1} \newcommand{\switchargs}[2]{#2#1}

\end{document}